\documentclass{amsproc}

\newtheorem{theorem}{Theorem}[section]
\newtheorem{lemma}[theorem]{Lemma}

\theoremstyle{definition}

\newtheorem{proposition}[theorem]{Proposition}

\theoremstyle{remark}

\numberwithin{equation}{section}

    \def\ds{\displaystyle}
    \DeclareMathOperator*{\diag}{diag}
    \DeclareMathOperator*{\sgn}{sgn}
    \DeclareMathOperator*{\supp}{supp}
    \DeclareMathOperator*{\Tr}{Tr}

\begin{document}

\title{Multiple orthogonal polynomial ensembles}
\author{Arno~B.~J. Kuijlaars}
\address{Department of Mathematics, Katholieke Universiteit
Leuven, Celestijnenlaan 200B, 3001 Leuven, Belgium}
\email{arno.kuijlaars@wis.kuleuven.be}

\thanks{The author was supported in part by FWO-Flanders project G.0427.09,
by K.U. Leuven research grant OT/08/33, by the Belgian Interuniversity Attraction Pole P06/02,
by the European
Science Foundation Program MISGAM, and by grant
MTM2008-06689-C02-01 of the Spanish Ministry of Science and Innovation.}

\date{February 6, 2009}

\dedicatory{Dedicated to Guillermo L\'opez Lagomasino,
on the occasion of his 60th birthday}


\begin{abstract}
Multiple orthogonal polynomials are traditionally studied
because of their connections to number theory and approximation theory.
In recent years they were found to be connected to certain
models in random matrix theory. In this paper we introduce the notion
of a multiple orthogonal polynomial ensemble (MOP ensemble) and derive some
of their basic properties. It is shown that Angelesco
and Nikishin systems give rise to MOP ensembles and that the
equilibrium problems that are associated with these systems have
a natural interpretation in the context of MOP ensembles.
\end{abstract}

\maketitle

\section{Introduction}

Multiple Orthogonal Polynomials (MOPs) were introduced and studied
for problems in analytic number theory (irrationality and transcendence
proofs). Later they appeared in approximation theory, most notably
in the theory of Hermite-Pad\'e approximation and in this context
they are also called Hermite-Pad\'e polynomials \cite{Apt1,Apt3,BL,DS1,DS2,GR1,GR2,GRS,Mah,Nut}.
MOPs were also studied from the point of view of new special functions
\cite{ABV,Borw,CV1,KVW,Sta,VC}. See the books \cite{Ism,NS} and the
survey papers \cite{Apt2,AS,VAs1,VAs2}
for these aspects of MOPs. Further developments in these directions
are reported in e.g.\ \cite{AKV,BCV,BK2,CCV,CV2,CV3,FIL,FLLS,LL,MV}.

Recently MOPs also appeared in a natural way in probability
theory and mathematical physics in certain models coming from
random matrix theory and non-intersecting paths. The connection
was first observed in \cite{BK0} where MOPs were used in
in a random matrix model with external source. In the Gaussian
case, the external source model has an equivalent interpretation
in terms of non-intersecting Brownian motions. The external
source model was further analyzed with the use of
multiple Hermite and multiple Laguerre polynomials in
\cite{ABK,Bai,BK1,BK3,DF1,IS,Ora,MM,Mo1,Mo2,Wan}, see also \cite{ABV,Bai,BK2,Des,LW}.
A related non-intersecting path model was studied in
\cite{KMW} using MOPs for modified Bessel weights that were introduced
earlier in \cite{CV1}.
The biorthogonal polynomials arising in the two matrix model were identified
as MOPs in \cite{KMc}.
For a special case they were asymptotically analyzed in \cite{DuK2,Mo3}.
The Cauchy two matrix model and their associated
Cauchy biorthogonal polynomials have a number of similar features
\cite{BGS1,BGS2}. MOPs were generelized to MOPs of mixed type in
\cite{AVV,DaK2,DKV,DeK}.

Asymptotic results were mainly obtained from an analysis of the Riemann-Hilbert
problem for MOPs, formulated by Van Assche et al.\ \cite{VGK}
as an extension of the Riemann-Hilbert problem for orthogonal
polynomials \cite{FIK}. The application of the Deift/Zhou steepest
descent analysis \cite{Dei} to the Riemann-Hilbert problem for MOPs
presents several interesting new features that however we will not discuss here.

It is the aim of this paper to give an introductory
account of MOPs from the point of view of determinantal
point processes. After discussing the definition and
some of the basic properties of MOPs we discuss a multiple
integral representatons for the type II MOPs, which is
essentially taken from \cite{BK0}.
Under a suitable
constant sign condition the formula can be interpreted as the
expectation value of the
random polynomial $\prod_{j=1}(z-x_j)$ with roots $x_1, \ldots, x_n$
from a determinantal point process (called a MOP ensemble)
on the real line.

The constant sign condition holds in particular for Angelesco and Nikishin systems.
For both of these systems we show that the joint p.d.f.\ of the
associated MOP ensemble takes on a particular nice form.
In the large $n$ limit it allows for a natural probabilistic interpretation
of the vector equilibrium problems
that are associated with Angelesco and Nikishin systems.

\section{Multiple orthogonal polynomials}
\subsection{Definitions}
Given weight functions $w_1, \ldots, w_p$ on $\mathbb R$ and a
multi-index $\vec{n} = (n_1, \ldots, n_p) \in \mathbb N^p$,
the type II MOP is a monic polynomial $P_{\vec{n}}$ of
degree $|\vec{n}| = n_1 + \cdots + n_p$ such that
\begin{equation} \label{eq:MOP2orthogonality}
    \int_{-\infty}^{\infty} P_{\vec{n}}(x) x^k w_j(x) dx = 0,
        \qquad k=0, \ldots, n_j-1, \quad j=1, \ldots, p.
        \end{equation}
Throughout we will write
\[ n = |\vec{n}| = n_1 + \cdots + n_p. \]

The conditions \eqref{eq:MOP2orthogonality} give a system of
$n$ linear equations for the $n$ free
coefficients of the polynomial $P_{\vec{n}}$ (recall that $P_{\vec{n}}$
is monic). If the system has a unique solution we say that
the multi-index $\vec{n}$ is normal (with respect to the
weights $w_1, \ldots, w_p$).

In this paper we mainly deal with the type II MOP, but at times
it is useful to consider the dual notion of type I MOPs as well.
These are polynomials $A_{\vec{n}}^{(j)}$, $j=1, \ldots, p$,
of degrees $\deg A_{\vec{n}}^{(j)} = n_j-1$, such that
the linear form
\begin{equation} \label{eq:MOP1linearform}
     Q_{\vec{n}}(x) = \sum_{j=1}^p A_{\vec{n}}^{(j)}(x) w_j(x)
     \end{equation}
satisfies
\begin{equation} \label{eq:MOP1orthogonality}
    \int_{-\infty}^{\infty} x^k Q_{\vec{n}}(x) dx = 0,
        \qquad k=0,1, \ldots, n-2.
        \end{equation}
If we supplement this with the normalizing condition
\begin{equation} \label{eq:MOP1normalization}
    \int_{-\infty}^{\infty} x^{n-1} Q_{\vec{n}}(x) dx = 1,
\end{equation}
then again we have a system of $n = |\vec{n}|$ linear equations for the in total
$n$ coefficients of the polynomials $A_{\vec{n}}^{(j)}$,
$j=1, \ldots, p$.

\subsection{Determinantal expressions}
Let
\[ c^{(j)}_k = \int_{-\infty}^{\infty} x^k w_j(x) dx \]
denote the $k$th moment of the weight $w_j$, and let
\[ H^{(j)}_{m,n} = \left( c^{(j)}_{k+l}\right)_{k=0, \ldots, m, l = 0, \ldots, n} \]
be the $(m+1) \times (n+1)$ Hankel matrix with the moments of $w_j$.
The conditions \eqref{eq:MOP1orthogonality} and \eqref{eq:MOP1normalization}
give rise to a linear system whose matrix has the block Hankel structure
\begin{equation} \label{eq:blockHankel}
    M_{\vec{n}} = \begin{bmatrix}
    H^{(1)}_{n-1, n_1-1} &
    H^{(2)}_{n-1, n_2-1} &
    \cdots &
    H^{(p)}_{n-1, n_p-1}
    \end{bmatrix}.
\end{equation}
Therefore the type I MOPs uniquely exist if and only if
\begin{equation} \label{eq:blockHankeldet}
    D_{\vec{n}} := \det M_{\vec{n}} =  \begin{vmatrix}
    H^{(1)}_{n-1, n_1-1} &
    H^{(2)}_{n-1, n_2-1} &
    \cdots &
    H^{(p)}_{n-1, n_p-1}
\end{vmatrix} \neq 0.
\end{equation}

The linear system arising from the type II conditions
\eqref{eq:MOP2orthogonality} has a matrix which is the transpose
of \eqref{eq:blockHankel}. Therefore the non-vanishing
of the determinant \eqref{eq:blockHankeldet} also
guarantees the existence and uniqueness of the type II MOP.

Suppose $D_{\vec{n}} \neq 0$. Then it is easy to see that
the type II MOP has the determinantal formula
\begin{equation} \label{eq:typeIIdeterminant}
    P_{\vec{n}}(x) = \frac{1}{D_{\vec{n}}}
\begin{vmatrix}
H^{(1)}_{n, n_1-1} &
H^{(2)}_{n,n_2-1} & \cdots &
 H^{(p)}_{n, n_p-1} &
 \begin{smallmatrix} \ds 1 \\[5pt] \ds x\\[5pt] \ds x^2 \\ \ds \vdots \\ \ds x^{n} \end{smallmatrix}
\end{vmatrix}.
\end{equation}
Indeed, the right-hand side of \eqref{eq:typeIIdeterminant} is a monic polynomial
of degree $n$. If we multiply the right-hand side of \eqref{eq:typeIIdeterminant}
by $x^k w_j(x)$ and integrate with respect to $x$,
we can perform these operations in the last column to obtain a determinant
with two equal columns if $k \leq n_j-1$. This proves the type II orthogonality
conditions \eqref{eq:MOP2orthogonality}.

The type I MOPs have a similar determinantal expression. For $j =1, \ldots, p$ we have
\begin{multline} \label{eq:typeIdeterminant}
    A_{\vec{n}}^{(j)}(x) =
    \frac{1}{D_{\vec{n}}} \times \\
        \begin{vmatrix}
        H^{(1)}_{n-2, n_1-1} &
        \cdots & H^{(j-1)}_{n-2,n_{j-1}-1} &
        H^{(j)}_{n-2, n_j-1} &
        H^{(j+1)}_{n-2,n_{j+1}-1} & \cdots &
        H^{(p)}_{n-2, n_p-1} \\[10pt]
        0 & \cdots & 0 &
        \begin{smallmatrix} \ds 1 & \ds x & \ds \cdots & \ds x^{n_j-1} \end{smallmatrix} &
        0 & \cdots & 0 \end{vmatrix}.
        \end{multline}

These and similar determinantal formulas have recently
been considered from the point of view of integrable systems in
\cite{AVV,Ber}.

\subsection{Multiple integral representation}

For what follows it is convenient to write
\[ N_j = \sum_{i=1}^j n_i, \qquad N_0 = 0 \]
and to introduce two sequences of functions $f_1, \ldots, f_{\vec{n}}$
and $g_1, \ldots, g_{\vec{n}}$ by
\begin{equation} \label{eq:deffj}
    f_j(x) = x^{j-1}, \qquad j = 1, \ldots, n
    \end{equation}
and
\begin{equation} \label{eq:defgj}
    g_{i + N_{j-1}}(x) =  x^{i-1} w_j(x), \qquad
    i = 1, \ldots, n_j, \quad  j= 1, \ldots, p.
    \end{equation}
Then the block Hankel matrix \eqref{eq:blockHankel} can be
written as
\begin{equation} \label{eq:defMn}
    M_{\vec{n}} = \begin{bmatrix} m_{j,k} \end{bmatrix}_{j,k=1, \ldots, n},
    \qquad m_{j,k} = \int_{-\infty}^{\infty} f_j(x) g_k(x) \, dx
    \end{equation}
and
\[ D_{\vec{n}} = \det M_{\vec{n}} =
    \det \begin{bmatrix} \int_{-\infty}^{\infty} f_j(x) g_k(x) \, dx \end{bmatrix}_{j,k=1, \ldots,\vec{n}}. \]

For general $m$ and $n = |\vec{n}|$ we also write
\begin{equation} \label{eq:defMmn}
    M_{m,n} = \begin{bmatrix} m_{j,k} \end{bmatrix}_{j=1, \ldots, m, k=1, \ldots, n},
    \end{equation}
so that we have by \eqref{eq:typeIIdeterminant}
\begin{equation} \label{eq:typeIIdeterminantM}
    P_{\vec{n}}(x) = \frac{1}{D_{\vec{n}}}
\begin{vmatrix} M_{n+1,n} &
 \begin{smallmatrix} \ds 1 \\[5pt] \ds x \\ \ds \vdots \\ \ds x^{n} \end{smallmatrix}
\end{vmatrix} \end{equation}
and by \eqref{eq:typeIdeterminant} and \eqref{eq:defgj}
\begin{equation} \label{eq:typeIdeterminantM}
    Q_{\vec{n}}(x) =
    \sum_{j=1}^p A_{\vec{n}}^{(j)}(x) w_j(x)
    = \frac{1}{D_{\vec{n}}}
        \begin{vmatrix}
        M_{n-1,n} \\
        \begin{smallmatrix} \ds g_1(x) & \ds g_2(x) & \ds \cdots & \ds g_{n}(x) \end{smallmatrix} \end{vmatrix}.
        \end{equation}

The following lemma is standard, see e.g.\ \cite[Proposition 2.10]{Joh}
where it is called a generalized Cauchy-Binet identity.
\begin{lemma}
We have
\begin{align} \label{eq:Dnasmultint}
    D_{\vec{n}} =  \frac{1}{n!} \int_{-\infty}^{\infty} \cdots \int_{-\infty}^{\infty}
        \det \begin{bmatrix} f_j(x_k) \end{bmatrix}_{j,k=1, \ldots, n}
    \cdot \det \begin{bmatrix} g_j(x_k) \end{bmatrix}_{j,k=1, \ldots, n}
        \, \prod_{k=1}^{n} dx_k.
    \end{align}
\end{lemma}
\begin{proof}
Expanding the two determinants on the right-hand side of \eqref{eq:Dnasmultint}
we get
\begin{align*}
    \det \begin{bmatrix} f_j(x_k) \end{bmatrix}_{j,k}
    \cdot \det \begin{bmatrix} g_j(x_k) \end{bmatrix}_{j,k}
     = \sum_{\sigma} \sum_{\tau} (-1)^{\sgn \sigma + \sgn \tau}
        \prod_{k=1}^{n} f_{\sigma(k)}(x_k) g_{\tau(k)}(x_k)
        \end{align*}
where the sums are for $\sigma$ and $\tau$ over the symmetric group $S_{n}$.
By \eqref{eq:defMn} the right-hand side of \eqref{eq:Dnasmultint} is equal to
\begin{equation} \label{eq:Dnasmultintproof1}
    \frac{1}{\vec{n}!} \sum_{\sigma} \sum_{\tau} (-1)^{\sgn \sigma + \sgn \tau}
    \prod_{k=1}^{n} m_{\sigma(k), \tau(k)}
    =  \frac{1}{n!} \sum_{\sigma} \sum_{\tau} (-1)^{\sgn (\sigma \circ \tau^{-1})}
    \prod_{k=1}^{n} m_{\sigma \circ \tau^{-1}(k), k}.
    \end{equation}
For any fixed $\sigma$, we have that $\sigma \circ \tau^{-1}$ runs through $S_{n}$
as $\tau$ runs through $S_{n}$. Hence
\begin{equation} \label{eq:Dnasmultintproof2}
    \sum_{\tau} (-1)^{\sgn (\sigma \circ \tau^{-1})}
    \prod_{k=1}^{\vec{n}} m_{\sigma \circ \tau^{-1}(k), k} = \det M_{\vec{n}} = D_{\vec{n}}.
    \end{equation}
The equality \eqref{eq:Dnasmultint} follows from \eqref{eq:Dnasmultintproof1}
and \eqref{eq:Dnasmultintproof2}.
\end{proof}

There is a similar multiple integral representation for the type II MOPs,
which was stated for a special case in \cite{BK0}, see also \cite{DF2}.
We emphasize that it is important here that $f_j(x) = x^{j-1}$.
\begin{proposition} \label{prop:Pnasmultint}
Assume $D_{\vec{n}} \neq 0$. Then the type II MOP has
the multiple integral representation
\begin{multline} \label{eq:Pnasmultint}
   P_{\vec{n}}(z) = \frac{1}{D_{\vec{n}} \cdot n!} \times \\
    \int_{-\infty}^{\infty} \cdots \int_{-\infty}^{\infty}
        \prod_{k=1}^{n} (z-x_k) \cdot \det \begin{bmatrix} f_j(x_k) \end{bmatrix}_{j,k=1, \ldots, n}
    \cdot \det \begin{bmatrix} g_j(x_k) \end{bmatrix}_{j,k=1, \ldots, n}
        \, \prod_{k=1}^{n} dx_k.
    \end{multline}
\end{proposition}
\begin{proof}
Since $f_j(x) = x^{j-1}$ we have that $\det \left[ f_j(x_k)\right]$ is
a Vandermonde determinant, and therefore
\[ \prod_{k=1}^{\vec{n}} (z-x_k) \cdot \det \begin{bmatrix} f_j(x_k) \end{bmatrix}_{j,k=1, \ldots, n}
    = \det \begin{bmatrix} f_j(x_k) \end{bmatrix}_{j,k=1, \ldots, n+1} \]
where we have put
\[ f_{n+1}(x) = x^{n}, \qquad \text{and} \qquad x_{n + 1} = z. \]
Thus, by expanding the determinant we have
\[ \prod_{k=1}^{n} (z-x_k) \cdot \det \begin{bmatrix} f_j(x_k) \end{bmatrix}_{j,k=1, \ldots, n}
    = \sum_{\sigma \in S_{n+1}} (-1)^{\sgn \sigma} \prod_{k=1}^{n}
        f_{\sigma(k)} (x_k) \cdot f_{\sigma(n+1)}(z) \]
and similarly
\[ \det \begin{bmatrix} g_j(x_k) \end{bmatrix}_{j,k=1, \ldots, n}
    = \sum_{\tau \in S_{n}} (-1)^{\sgn \tau} \prod_{k=1}^{n} g_{\tau(k)}(x_k). \]

Integrating the product of the two above expressions with respect to
$x_1, \ldots, x_{n}$ we obtain
\begin{multline} \label{eq:Pnasmultintproof1}
     \int_{-\infty}^{\infty} \cdots \int_{-\infty}^{\infty}
        \prod_{k=1}^{n} (z-x_k) \cdot \det \begin{bmatrix} f_j(x_k) \end{bmatrix}_{j,k=1, \ldots, n}
    \cdot \det \begin{bmatrix} g_j(x_k) \end{bmatrix}_{j,k=1, \ldots, n}
        \, \prod_{k=1}^{n} dx_k \\
        =
        \sum_{\sigma \in S_{n+1}} \sum_{\tau \in S_{n}}
        (-1)^{\sgn \sigma + \sgn \tau} \prod_{k=1}^{n} m_{\sigma(k), \tau(k)} \cdot f_{\sigma(n+1)}(z)
        \\
        =
        \sum_{\tau \in S_{n}} \sum_{\sigma \in S_{n+1}}
        (-1)^{\sgn (\sigma \circ \tau^{-1})}
        \prod_{k=1}^{n} m_{\sigma \circ \tau^{-1}(k), k} \cdot z^{\sigma(n+1)-1},
        \end{multline}
where we used the definition of $m_{j,k}$ as given in \eqref{eq:defMn} also for $j = n +1$.

For each fixed $\tau \in S_{n}$ we have that the sum over $\sigma$
in \eqref{eq:Pnasmultintproof1} is equal to the determinant in
the right-hand side of \eqref{eq:typeIIdeterminantM} and the proposition follows.
\end{proof}

In an analogous way we find the following multiple integral representation
for the linear form of type I MOPs, which is due to Desrosiers and Forrester \cite{DF2}.
\begin{proposition} \label{prop:Qnasmultint}
Assume $D_{\vec{n}} \neq 0$. Then the linear form of type I MOPs satisfies
\begin{multline} \label{eq:Qnasmultint}
    \int_{-\infty}^{\infty} \frac{Q_{\vec{n}}(x) w_j(x)}{z-x} dx
   = \frac{1}{D_{\vec{n}} \cdot n!} \times \\
    \int_{-\infty}^{\infty} \cdots \int_{-\infty}^{\infty}
        \prod_{k=1}^{n} (z-x_k)^{-1} \cdot \det \begin{bmatrix} f_j(x_k) \end{bmatrix}_{j,k=1, \ldots, n}
    \cdot \det \begin{bmatrix} g_j(x_k) \end{bmatrix}_{j,k=1, \ldots, n}
        \, \prod_{k=1}^{n} dx_k.
    \end{multline}
\end{proposition}
\begin{proof}
Here we use the property
\[ \prod_{k=1}^{n} (z-x_k)^{-1} \cdot \det \begin{bmatrix} f_j(x_k) \end{bmatrix}_{j,k=1, \ldots, n}
    = \det \begin{bmatrix} f_j(x_k) \end{bmatrix}_{j,k=1, \ldots, n+1} \]
where now we put
\[ f_{\vec{n} + 1}(x) = \frac{1}{z-x}, \qquad \text{and} \qquad x_{n + 1} = z. \]
The rest of the proof follows along the same lines
as the proof of Proposition \ref{prop:Pnasmultint}. We omit the details, see also \cite{DF2}.
\end{proof}

\section{MOP ensembles}

\subsection{Probabilistic interpretation}
The multiple integral representations \eqref{eq:Dnasmultint}, \eqref{eq:Pnasmultint}
and \eqref{eq:Qnasmultint} have a natural probabilistic interpretation
in case the product of determinants
\[ \det \begin{bmatrix} f_j(x_k) \end{bmatrix}_{j,k=1, \ldots, n}
    \cdot \det \begin{bmatrix} g_j(x_k) \end{bmatrix}_{j,k=1, \ldots, n} \]
is of a fixed sign for $(x_1, \ldots, x_{n}) \in \mathbb R^{n}$.
That is, if it is always $\geq 0$ or always~$\leq 0$.
Indeed, in that case it follows by \eqref{eq:Dnasmultint} that
\begin{equation} \label{eq:biorthogonalpdf}
    \mathcal P(x_1,\ldots, x_n) =
        \frac{1}{Z_n}
        \det \left[f_j(x_k)\right]_{j,k=1,\ldots,n} \cdot \det \left[g_j(x_k) \right]_{j,k=1,\ldots,n}
        \end{equation}
        is a probability density function on $\mathbb R^n$, where
\begin{equation} \label{eq:biorthogonalZn}
    Z_n = D_n n! \end{equation}
is the normalizing constant (also called partition function in statistical mechanics
literature), so
that $\int \cdots \int \mathcal P(x_1, \ldots, x_n) dx_1 \cdots dx_n = 1$.

The multiple integral representations \eqref{eq:Pnasmultint} and \eqref{eq:Qnasmultint} then show that
\begin{align} \label{eq:Pnasexpection}
    P_{\vec{n}}(z) & = \mathbb E \left[ \prod_{k=1}^n (z-x_k) \right], \qquad z \in \mathbb C, \\
    Q_{\vec{n}}(z) & = \mathbb E \left[ \prod_{k=1}^n (z-x_k)^{-1} \right], \qquad
    z \in \mathbb C \setminus \mathbb R, \label{eq:Qnasexpectation}
    \end{align}
where the mathematical expectation is taken with respect to the p.d.f.\ \eqref{eq:biorthogonalpdf}.

Thus $P_{\vec{n}}(z)$ is the average of the polynomials $\prod_{k=1}^n (z-x_k)$ where
the roots $x_1, \ldots, x_n$ are distributed according to \eqref{eq:biorthogonalpdf}.
In cases where the distribution \eqref{eq:biorthogonalpdf} can be interpreted as the
eigenvalue distribution of a random matrix ensembles, one would call $P_{\vec{n}}$
the average characteristic polynomial.

\subsection{Biorthogonal ensembles}

A biorthogonal ensemble, see \cite{Bor}, is a probability density function on $\mathbb R^n$ of the form
\eqref{eq:biorthogonalpdf} with certain given functions $f_1, \ldots, f_n$ and $g_1, \ldots, g_n$,
not necessarily of the form \eqref{eq:deffj} and \eqref{eq:defgj}.
The p.d.f.\ is invariant under permutations of variables. We think of
the ensemble as giving us $n$ random points or particles $x_j$ on the real line, and so
it is a random point process.

A biorthogonal ensemble is a special case of a determinantal point process,
see e.g.\ \cite{Joh, Sos},
This means that there is a correlation kernel $K_n(x,y)$ so that
    \[ \mathcal P(x_1, \ldots, x_n) =
        \frac{1}{n!} \det \left[K_n(x_j,x_k)\right]_{j,k=1, \ldots, n} \]
    and so that marginal densities ($m$ point correlation functions) are determinants
    \begin{align*}
        \underbrace{\int_{-\infty}^{\infty} \cdots \int_{-\infty}^{\infty}}_{n-m \textrm{ times}} \mathcal P(x_1, \ldots, x_n) dx_{m+1} \cdots dx_n
         =
        \frac{(n-m)!}{n!} \det \left[K_n(x_j,x_k)\right]_{j,k=1, \ldots, m}.
    \end{align*}
Taking for example $m=1$ we have that $\frac{1}{n} K_n(x,x)$
is the mean  density of points, that is
\[ \frac{1}{n} \int_{a}^b K_n(x,x) dx \]
is the expected fraction of points lying in the interval $[a,b]$.

In a biorthogonal ensemble, the correlation kernel can be written as a bordered determinant
\begin{equation} \label{eq:Knasdeterminant}
    K_n(x,y) =
        \frac{-1}{\det M_n} \begin{vmatrix}
            M_n & \begin{smallmatrix} \ds f_1(x) \\ \ds \vdots \\ \ds f_n(x) \end{smallmatrix} \\
            \begin{smallmatrix} \ds g_1(y) & \ds \cdots &  \ds g_n(y) \end{smallmatrix} & 0
            \end{vmatrix} \end{equation}
    where $M_n$ is the matrix
    \[ M_n = \begin{bmatrix} m_{j,k} \end{bmatrix}_{j,k=1, \ldots, n},
        \qquad m_{j,k} =  \int_{-\infty}^{\infty} f_j(x) g_k(x) \, dx  \]

In the formulation of the biorthogonal ensemble \eqref{eq:biorthogonalpdf},
we have some freedom in choosing the functions $f_1, \ldots, f_n$. and $g_1,\ldots, g_n$.
Indeed, if $\phi_1, \ldots, \phi_n$ and $\psi_1,\ldots, \psi_n$ are functions
with the same linear span as the $f_j$'s and $g_j$'s, respectively,
then we could use these functions instead. A particular nice form
appears if the functions $\phi_j$ and $\psi_k$ are biorthogonal, i.e.,
\[ \int_{-\infty}^{\infty} \phi_j(x) \psi_k(x) dx = \delta_{j,k}. \]
Then the representation \eqref{eq:Knasdeterminant} reduces to
\begin{align} \label{eq:Knasdeterminant2} K_n(x,y) & =
        - \begin{vmatrix}
            I_n & \begin{smallmatrix} \ds \phi_1(x) \\ \ds \vdots \\ \ds \phi_n(x) \end{smallmatrix} \\
            \begin{smallmatrix} \ds \psi_1(y) & \ds \cdots & \ds \psi_n(y) \end{smallmatrix} & 0
            \end{vmatrix}  = \sum_{j=1}^n \phi_j(x) \psi_j(y). \end{align}

\subsection{OP ensembles}
If $f_j(x) = g_j(x) = x^{j-1} \sqrt{w(x)}$, $j=1, \ldots, n$ for some
non-negative weight function $w$ on $\mathbb R$, then
\begin{equation} \label{eq:OPensemblepdf}
    \frac{1}{Z_n} \det \left[f_j(x_k)\right] \cdot \det \left[g_j(x_k) \right]
        = \frac{1}{Z_n} \prod_{1 \leq j < k \leq n} (x_k-x_j)^2 \cdot \prod_{k=1}^n w(x_k)
        \end{equation}
which is indeed of constant sign.
This is the form of the joint p.d.f.\ for the eigenvalues of a unitary random matrix
ensemble
\[ \frac{1}{\tilde{Z}_n} \exp( - \Tr V(H)) dH, \qquad w(x) = e^{-V(x)}, \]
defined on $n \times n$ Hermitian matrices $H$, see \cite{Dei}.

In this case the biorthogonal functions take the form
\begin{equation} \label{eq:OPpolynomials}
    \phi_j(x) = \psi_j(x) = p_{j-1}(x) \sqrt{w(x)}
    \end{equation}
where  $p_{j-1}$ is the orthonormal polynomial of degree $j-1$ with
respect to the weight $w$ on $\mathbb R$, and by \eqref{eq:Knasdeterminant2},
\[ K_n(x,y) = \sqrt{w(x)} \sqrt{w(y)}
        \sum_{j=0}^{n-1} p_j(x) p_j(y) \]
        is the correlation kernel, which in this situation is
        also called the OP kernel or the Christoffel-Darboux kernel.

The OPs are characterized by a $2 \times 2$ matrix valued Riemann-Hilbert problem
due to Fokas, Its, and Kitaev \cite{FIK},
     \begin{itemize}
    \item $Y : \mathbb C \setminus \mathbb R \to \mathbb C^{2 \times 2}$ is analytic,
    \item $Y_+(x) = Y_-(x) \begin{pmatrix} 1 & w(x) \\ 0 & 1 \end{pmatrix}$,
        \quad $x \in \mathbb R$,
    \item $Y(z) = \left(I_{2} + \mathcal{O}(1/z) \right)
        \diag \begin{pmatrix} z^{n} & z^{-n} \end{pmatrix}$
        \quad as $z \to \infty$.
        \end{itemize}
The correlation kernel for the OP ensemble can be given directly in terms of
the solution $Y$ of the RH problem
\begin{equation} \label{eq:KninRHP}
    K_n(x,y) = \frac{1}{2\pi i(x-y)} \sqrt{w(x)} \sqrt{w(y)}
        \begin{pmatrix} 0 & 1 \end{pmatrix}
        Y_+(y)^{-1} Y_+(x) \begin{pmatrix} 1 \\ 0 \end{pmatrix}.
        \end{equation}
This follows from an explicit formula for $Y$ in terms of
the orthogonal polynomials $p_n$ and $p_{n-1}$,
and the Christoffel-Darboux formula  for orthogonal polynomials.

\subsection{MOP ensembles}
We have a MOP ensemble if $f_1, \ldots, f_n$ and $g_1, \ldots, g_n$
are given by  \eqref{eq:deffj} and \eqref{eq:defgj}, and
if
\begin{equation} \label{eq:MOPensemblecondition}
    \det \begin{bmatrix} f_j(x_k) \end{bmatrix}_{j,k=1, \ldots, n}
    \cdot \det \begin{bmatrix} g_j(x_k) \end{bmatrix}_{j,k=1, \ldots, n}
    \end{equation}
has constant sign. The case $p=1$ reduces essentially to the OP case.

For a MOP ensemble we have that the correlation kernel $K_n$ given by
the determinant \eqref{eq:Knasdeterminant} has another expression in terms of the RH problem for
multiple orthogonal polynomials. MOPs (with $p$ weights) satisfy a
$(p+1) \times (p+1)$ matrix valued RH problem \cite{VGK}
     \begin{itemize}
    \item $Y : \mathbb C \setminus \mathbb R \to \mathbb C^{(p+1)\times (p+1)}$ is analytic,
    \item $Y_+(x) = Y_-(x) \begin{pmatrix} 1 & w_1 & w_2 & \cdots & w_p \\
    0 & 1 & 0 &  & 0 \\
     &  0 & \ddots & & \vdots \\
     &  &  & \ddots & 0 \\
    0 &  & & 0 & 1  \end{pmatrix}$,
        \quad $x \in \mathbb R$,
    \item $Y(z) = \left(I_{p+1} + \mathcal{O}(1/z) \right)
        \diag \begin{pmatrix} z^{n}  &
        z^{-n_1} & \cdots & z^{-n_p} \end{pmatrix}$  as $z \to \infty$.
        \end{itemize}
The correlation kernel for the MOP ensemble is given as follows in terms of
    the solution $Y$ of the RH problem
\begin{align} \label{eq:KninRHPforMOP}
    K_n(x,y) = \frac{1}{2\pi i(x-y)}
        \begin{pmatrix} 0 & w_{1}(y) & \cdots & w_{p}(y) \end{pmatrix}
        Y_+(y)^{-1}   Y_+(x) \begin{pmatrix} 1 \\ 0 \\ \vdots \\ 0 \end{pmatrix}.
        \end{align}
The proof is based on a Christoffel-Darboux formula for MOPs,
see \cite{BK0} for the case $p=2$ and \cite{DaK1} for general $p$.
An extension to MOP of mixed type is given in \cite{DaK2}.

\section{Special classes}

The condition that \eqref{eq:MOPensemblecondition} is of constant sign
is automatically satisfied in the OP case, but it becomes relevant in the MOP case.
We call it the MOP ensemble condition.
It is of interest to identify classes for which the MOP ensemble condition holds.
In the literature on Hermite-Pad\'e approximation a number of special
classes of MOPs were identified for which rather complete convergence results
could be established. These are in particular the Angelesco systems \cite{Apt1,GR1}
and the Nikishin systems, see e.g.\ \cite{Apt3,BL,DS2,GRS,LL}.
It turns out that for these special classes
the MOP ensemble condition holds. Before we turn to that, we
make some preliminary observations.

\subsection{Preliminary observations}
The first observation is that the product of determinants
\eqref{eq:MOPensemblecondition} is invariant under permutations of the
$x_k$'s. It is also clear that \eqref{eq:MOPensemblecondition} is zero
in case two or more of the $x_k$'s coincide. Therefore we may restrict ourselves
to strictly ordered sets of points
\begin{equation} \label{eq:xkordered}
    x_1 < x_2 < \cdots < x_n.
    \end{equation}
The second observation is that the first factor in \eqref{eq:MOPensemblecondition}
is a Vandermonde determinant (due to the fact that $f_j(x) = x^{j-1}$)
\[  \det \begin{bmatrix} f_j(x_k) \end{bmatrix}_{j,k=1, \ldots, n}
    = \prod_{j < k} (x_k -x_j) \]
which is positive for ordered points \eqref{eq:xkordered}.
Therefore the MOP ensemble condition comes down to the
condition stated in the following lemma.
\begin{lemma} \label{lem:MOPcondition}
The MOP ensemble condition is satisfied if and only if
either
\begin{equation} \label{eq:detgjpositive}
    \det \begin{bmatrix} g_j(x_k) \end{bmatrix}_{j,k=1, \ldots, n} \geq 0
    \end{equation}
whenever $x_1 < x_2 < \cdots < x_n$, or
\begin{equation} \label{eq:detgjnegative}
    \det \begin{bmatrix} g_j(x_k) \end{bmatrix}_{j,k=1, \ldots, n} \leq 0
    \end{equation}
whenever $x_1 < x_2 < \cdots < x_n$.
\end{lemma}

Since Vandermonde-like products will appear frequently in what follows
we use the abbreviations
\begin{equation} \label{eq:DeltaX}
    \Delta(X) = \prod_{1 \leq j < k \leq n} (x_k -x_j)
    \end{equation}
if $X = (x_1, \ldots, x_n)$, and
\begin{equation} \label{eq:DeltaXY}
    \Delta(X,Y) = \prod_{k=1}^n \prod_{j=1}^m (x_k - y_j)
    \end{equation}
if $X = (x_1, \ldots, x_n)$ and $Y = (y_1, \ldots, y_m)$.

\subsection{Angelesco ensemble}

The weights $w_1, \ldots, w_p$ are an Angelesco system
if there are disjoint intervals $\Gamma_1, \ldots, \Gamma_p \subset \mathbb R$,
such that
\[ \supp(w_j) \subset \Gamma_j, \qquad j=1, \ldots, p. \]
We write $\Gamma_j = [\alpha_j, \beta_j]$ and without loss of generality
we assume that
\begin{equation} \label{eq:Deltaordered}
    \beta_j < \alpha_{j+1}, \qquad \text{for } j=1, \ldots, p-1.
    \end{equation}
We extend $w_{j}$ to all of $\mathbb R$ by defining
    $w_{j}(x) = 0$ if $x \in \mathbb R \setminus \Gamma_j$.

An Angelesco system always gives rise to a MOP ensemble.
Indeed it is easy to see that in the  Angelesco case $\det \begin{bmatrix} g_j(x_k) \end{bmatrix}$
is of block form, and it can only be non-zero if $n_j$ of the points $x_k$
belong to $\Gamma_j$ for $j=1, \ldots, p$, and so this is what
we will assume. Then writing $N_j = \sum_{i=1}^j n_i$, $N_0 = 0$,
we have that
\[ x_k^{(j)} := x_{N_{j-1}+k} \in \Gamma_j, \qquad  k =1, \ldots, n_j, \quad j=1, \ldots, p. \]
Because of the orderings \eqref{eq:xkordered} and \eqref{eq:Deltaordered}
the determinant
then has a block diagonal form where the $i$th block is
\[ \begin{pmatrix}  w_i(x_1^{(i)}) & \cdots &  \cdots & w_i(x_{n_i}^{(i)}) \\
    x_{1}^{(i)} w_i(x_1^{(i)}) & \cdots &  \cdots & x_{n_i}^{(i)} w_i(x_{n_i}^{(i)}) \\
    \vdots &  & & \vdots \\
    \left(x_{1}^{(i)}\right)^{n_i-1} w_i(x_1^{(i)}) & \cdots &  \cdots &
    \left(x_{n_i}^{(i)}\right)^{n_i-1} w_i(x_{n_i}^{(i)})
    \end{pmatrix} \]
whose determinant is
\[ \prod_{1 \leq j < k \leq n_i} (x_k^{(i)}-x_j^{(i)})  \cdot \prod_{k=1}^{n_i} w_i(x_k^{(i)})
    = \Delta(X^{(i)}) \prod_{k=1}^{n_i} w_i(x_k^{(i)}) \]
where $X^{(i)} = (x_1^{(i)}, \ldots, x_{n_i}^{(i)})$.
The result is that
\begin{align}
    \det \begin{bmatrix} g_j(x_k) \end{bmatrix} & = \prod_{i=1}^p \left(
        \Delta(X^{(i)}) \cdot \prod_{k=1}^{n_i} w_i(x_k^{(i)}) \right)
            \end{align}
and this is $\geq 0$ for every choice of $x_1 < \cdots < x_n$.

Thus, by Lemma \eqref{lem:MOPcondition}, an Angelesco system gives rise to a MOP ensemble and we call it
an Angelesco ensemble. We see from the above calculation that the joint p.d.f.\
in an Angelesco ensemble is
\begin{multline} \label{eq:Angelescopdf}
    \frac{1}{Z_n} \det [ f_j(x_k)]  \det[g_j(x_k)] \\
      = \frac{1}{Z_n}
    \prod_{i=1}^p  \Delta(X^{(i)})^2 \cdot
    \prod_{1 \leq i < j \leq p} \Delta(X^{(i)}, X^{(j)}) \cdot
    \prod_{i=1}^p \prod_{k=1}^{n_i} w_i(x_k^{(i)}).
        \end{multline}

\subsection{AT ensemble}
Assume that $w_1, \ldots, w_p$ are weights defined on a
fixed interval $\Gamma \subset \mathbb R$.
Then $w_1, \ldots, w_p$ are an AT system on $\Gamma$
if the functions $g_j$ are an algebraic Chebyshev system on $\Gamma$.
This means that every non-trivial linear combination
\[ \sum_{j=1}^n \lambda_j g_j \]
has at most $n-1$ zeros in $\Gamma$.
Equivalently, an algebraic Chebyshev system means that
\begin{equation} \label{eq:ATdeterminant}
    \det \begin{bmatrix} g_j(x_k) \end{bmatrix} \neq 0
    \end{equation}
for every choice of distinct points $x_k$ in $\Gamma$.
The property of being an AT system also depends on the multi-indices
$n_1, \ldots, n_p$.

If the weights $w_j$ in an AT system are continuous functions on $\Gamma$,
then it clearly follows from \eqref{eq:ATdeterminant} by continuity
that $\det \begin{bmatrix} g_j(x_k) \end{bmatrix}$
has constant sign (either $> 0$ or $< 0$) whenever the $x_k$
are strictly ordered points in $\Gamma$.
Therefore, in that case, we have a MOP ensemble by Lemma \ref{lem:MOPcondition},
which we
we will call an AT ensemble.

\subsection{Nikishin ensemble}

\subsubsection{Definition of a Nikishin system}
Certain AT systems with special properties were first described by Nikishin \cite{Nik}
and are therefore called Nikishin systems. We state it first for $p=2$
continuous weight functions $w_1, w_2$ defined on an interval
$\Gamma_1 \subset \mathbb R$.

The Nikishin assumption is that the ratio $w_2/w_1$
can be written as a Markov function for a non-negative weight function (or more generally
a measure) supported on an interval $\Gamma_2$, disjoint from $\Gamma_1$,
that is, if
 \begin{equation} \label{eq:Nikishinratio1}
    \frac{w_2(x)}{w_1(x)} = \pm \int_{\Gamma_2} \frac{v(y)}{x-y}  dy,
        \qquad x \in \Gamma_1,
    \end{equation}
where $v(s)$ is a non-negative weight function with
\[ \supp(v) = \Gamma_2, \qquad \Gamma_2 \cap \Gamma_1 = \emptyset. \]
We choose the $+$ sign in \eqref{eq:Nikishinratio1}
if $\Gamma_2$ lies to the left of $\Gamma_1$; otherwise we choose the $-$ sign.
Then we call $w_1, w_2$ a Nikishin system on $\Gamma_1$
for the intervals $\Gamma_1, \Gamma_2$.

A Nikishin system with $p \geq 3$ weights is defined inductively.
Suppose $\supp(w_j) = \Gamma_1$ for all $j =1, \ldots, p$, where $\Gamma_1$ is an interval. Suppose
\begin{equation} \label{eq:Nikishinratio2}
    \frac{w_{j}(x)}{w_1(x)} = \pm \int_{\Gamma_2} \frac{v_j(y)}{x-y}  dy,
    \qquad x \in \Gamma_1,  \qquad j=2, \ldots, p,
        \end{equation}
where $\Gamma_2 \cap \Gamma_1 = \emptyset$ and where
$v_2, \ldots, v_p$ is a Nikishin system on $\Gamma_2$ for
the intervals $\Gamma_2, \ldots, \Gamma_p$.
Then we call $w_1, \ldots, w_p$ a Nikishin system on $\Gamma_1$ for the intervals
$\Gamma_1, \ldots, \Gamma_p$.

Note that in a Nikishin system two consecutive intervals $\Gamma_j$ and $\Gamma_{j+1}$
are disjoint. However, if $|j-k| \geq 2$, then $\Gamma_j$ and $\Gamma_k$  may
very well have a non-empty intersection.

This construction might not seem very natural at first sight, but it is actually
a very beautiful structure. A main result is that for multi-indices
$\vec{n} = (n_1, \ldots, n_p)$ such that
\[ n_j \geq n_{j+1} -1, \qquad j =1, \ldots, p-1 \]
a Nikishin system is an AT system, see \cite{NS},
and therefore the type I and type II
MOPs exist. As we have seen in the previous subsection, there is
also an associated MOP ensemble, which we call a Nikishin ensemble.

\subsubsection{Nikishin ensemble with $2$ weights}

Here we show that a Nikishin ensemble has a natural interpretation
as the marginal distribution of an extended ensemble.
The following calculations are due to Coussement and Van Assche \cite{CV1}.

For reasons of clarity we take $p=2$ and we assume that $\Gamma_2$ is
to the left of $\Gamma_1$.
Then for $x_1, \ldots, x_n$ in $\Gamma_1$ we have
\begin{align*} \det \begin{bmatrix} g_j(x_k)\end{bmatrix}
    & = \begin{vmatrix} w_1(x_1) & w_1(x_2) & \cdots & w_1(x_n) \\
        x_1 w_1(x_1) &  & &  x_n w_1(x_n) \\
        \vdots &  & & \vdots \\
        x_1^{n_1-1} w_1(x_1) & \cdots & \cdots & x_n^{n_1-1} w_1(x_n) \\
        w_2(x_1) & w_2(x_2) & \cdots & w_2(x_n) \\
        \vdots & & & \vdots \\
        x_1^{n_2-1} w_2(x_1) &  \cdots & \cdots & x_n^{n_2-1} w_2(x_n)
        \end{vmatrix} \\
    &= \prod_{k=1}^n w_1(x_k)  \begin{vmatrix} 1 & 1 & \cdots & 1 \\
        x_1 &  & &  x_n  \\
        \vdots &  & & \vdots \\
        x_1^{n_1-1} & \cdots & \cdots & x_n^{n_1-1}  \\
        \frac{w_2(x_1)}{w_1(x_1)} & \frac{w_2(x_2)}{w_1(x_2)} & \cdots & \frac{w_2(x_n)}{w_1(x_n)} \\
        \vdots & & & \vdots \\
        x_1^{n_2-1} \frac{w_2(x_1)}{w_1(x_1)} &  \cdots & \cdots & x_n^{n_2-1} \frac{w_2(x_n)}{w_1(x_n)}
        \end{vmatrix}
        \end{align*}
Now we replace each ratio $\frac{w_2(x_k)}{w_1(x_k)}$ by the integral
\eqref{eq:Nikishinratio1}, we use $y_j$ as the integration variable
in row $n_1 + j$, and we take the integrals as well as the factors $v(y_j)$
out of the determinant, to obtain
\[ \prod_{k=1}^n w_1(x_k)
    \underbrace{\int_{\Gamma_2} \cdots \int_{\Gamma_2}}_{n_2 \text{ times}}
        \prod_{j=1}^{n_2} v(y_j) \begin{vmatrix} 1 & 1 & \cdots & 1 \\
        \vdots &  & & \vdots \\
        x_1^{n_1-1} & \cdots & \cdots & x_n^{n_1-1}  \\
        \frac{1}{x_1-y_1}  & \frac{1}{x_2-y_1}  & \cdots & \frac{1}{x_n-y_1}  \\
        \frac{x_1}{x_1-y_2} & &  & \frac{x_n}{x_n-y_2} \\
        \vdots & & & \vdots \\
        \frac{x_1^{n_2-1}}{x_1-y_{n_2}} &  \cdots & \cdots & \frac{x_n^{n_2-1}}{x_n-y_{n_2}}
        \end{vmatrix}
            \prod_{j=1}^{n_2} d y_j. \]
Since
\[ n_1 \geq n_2 -1 \]
we can perform elementary row operations to reduce the remaining determinant
to
\[ \begin{vmatrix} 1 & 1 & \cdots & 1 \\
        \vdots &  & & \vdots \\
        x_1^{n_1-1} & \cdots & \cdots & x_n^{n_1-1}  \\
        \frac{1}{x_1-y_1}  & \frac{1}{x_2-y_1}  & \cdots & \frac{1}{x_n-y_1}  \\
        \frac{y_2}{x_1-y_2} & &  & \frac{y_2}{x_n-y_2} \\
        \vdots & & & \vdots \\
        \frac{y_2^{n_2-1}}{x_1-y_{n_2}} &  \cdots & \cdots & \frac{y_2^{n_2-1}}{x_n-y_{n_2}}
        \end{vmatrix}
        = \prod_{j=1}^{n_2} y_j^{j-1}
    \begin{vmatrix} 1 & 1 & \cdots & 1 \\
        \vdots &  & & \vdots \\
        x_1^{n_1-1} & \cdots & \cdots & x_n^{n_1-1}  \\
        \frac{1}{x_1-y_1}  & \frac{1}{x_2-y_1}  & \cdots & \frac{1}{x_n-y_1}  \\
        \frac{1}{x_1-y_2} & &  & \frac{1}{x_n-y_2} \\
        \vdots & & & \vdots \\
        \frac{1}{x_1-y_{n_2}} &  \cdots & \cdots & \frac{1}{x_n-y_{n_2}}
        \end{vmatrix}
        \]
which is a mixture of a Vandermonde and a Cauchy determinant. It can be evaluated
to give
\[  \prod_{j=1}^{n_2} y_j^{j-1} \cdot \frac{\Delta(X) \Delta(Y)}{\Delta(X,Y)}. \]
where $X = (x_1,\ldots, x_n)$ and $Y = (y_1, \ldots, y_{n_2})$.
Thus
\begin{multline} \label{eq:Nikishindetg2}
    \det \begin{bmatrix} g_j(x_k)\end{bmatrix} \\
    = \prod_{k=1}^n w_1(x_k)  \Delta(X)
    \int_{\Gamma_2} \cdots \int_{\Gamma_2}
        \prod_{j=1}^{n_2} v(y_j) \prod_{j=1}^{n_2} y_j^{j-1}
         \frac{\Delta(Y)}{\Delta(X,Y)}
            \prod_{j=1}^{n_2} d y_j. \end{multline}
Now we symmetrize the multiple integral with respect to the integration variables $y_j$
(which is a standard trick in determinantal point processes).
That is, for any permutation
$\sigma \in S_{n_2}$ we make the change of variables $y_j \mapsto y_{\sigma(j)}$ and we average over all
permutations $\sigma$ in $S_{n_2}$. Using the fact that
\[ \sum_{\sigma \in S_{n_2}} (-1)^{\sgn \sigma} \prod_{j=1}^{n_2} y_{\sigma_j}^{j-1}
    = \det\begin{bmatrix} y_k^{j-1} \end{bmatrix}_{j,k=1, \ldots, n_2}
    = \Delta(Y) \]
we then obtain that  \eqref{eq:Nikishindetg2} is equal to
\begin{multline} \label{eq:Nikishindetg3}
    \det \begin{bmatrix} g_j(x_k)\end{bmatrix} \\
    =  \frac{1}{n_2!} \prod_{k=1}^n w_1(x_k) \Delta(X)
    \int_{\Gamma_2} \cdots \int_{\Gamma_2}
        \prod_{j=1}^{n_2} v(y_j)
         \frac{\Delta(Y)^2}{\Delta(X,Y)}
            \prod_{j=1}^{n_2} d y_j. \end{multline}

The joint p.d.f.\ for the Nikishin ensemble is therefore
(since $\det \begin{bmatrix} f_j(x_k)\end{bmatrix} = \Delta(X)$)
\begin{multline} \label{eq:NikishinPDF}
    \frac{1}{Z_n} \det \begin{bmatrix} f_j(x_k)\end{bmatrix}
    \det \begin{bmatrix} g_j(x_k)\end{bmatrix} \\
    = \frac{1}{Z_n n_2!} \prod_{k=1}^n w_1(x_k)  \Delta(X)^2
    \int_{\Gamma_2} \cdots \int_{\Gamma_2}
        \prod_{j=1}^{n_2} v(y_j)
         \frac{\Delta(Y)^2}{\Delta(X,Y)}
            \prod_{j=1}^{n_2} d y_j. \end{multline}

By dropping the integrals over the $y_j$ variables, we can view
\eqref{eq:NikishinPDF} as
a marginal density of an extended ensemble defined by the joint p.d.f.\
\begin{multline} \label{eq:NikishinPDFextended}
    \mathcal P_{ext}(x_1, \ldots, x_n, y_1, \ldots, y_{n_2}) \\
    = \frac{1}{Z_n n_2!} \prod_{k=1}^n w_1(x_k)
        \prod_{j=1}^{n_2} v(y_j) \cdot
         \frac{\Delta(X)^2 \cdot \Delta(Y)^2}{\Delta(X,Y)}
         \end{multline}
defined for $x_1, \ldots, x_n \in~\Gamma_1$ and $y_1, \ldots, y_{n_2} \in \Gamma_2$.
Note that the factor $\Delta(X,Y)$ in  \eqref{eq:NikishinPDFextended}
is positive, since $\Gamma_2$ lies to the left of $\Gamma_1$ so
that $x_k > y_j$ for every $k = 1, \ldots, n$
and $j=1, \ldots, n_2$.

\subsubsection{Nikishin ensemble with $p \geq 2$ weights}

The above considerations can be extended to general $p \geq 2$.
Let $w_1, \ldots, w_p$ be a Nikishin system with $p$ weights
for the intervals $\Gamma_1, \ldots, \Gamma_p$. Assume that
\[ n_j \geq n_{j+1} - 1, \qquad \text{for } j=1, \ldots, p-1. \]
Then the joint p.d.f.\ for the Nikishin ensemble is a marginal
density of an extended ensemble defined on
\[ \Gamma_1^{N_1} \times \Gamma_2^{N_2} \times \cdots \times \Gamma_p^{N_p},
    \qquad \text{where } N_j = \sum_{i=j}^p n_i, \]
with joint p.d.f.\ of the form
\begin{align} \label{eq:NikishinExtendedPDF}
\frac{1}{\tilde{Z}_n}
    \prod_{j=1}^p \prod_{k=1}^{N_j} w^{(j)}(x_k^{(j)})
         \cdot
         \frac{\prod_{j=1}^p \Delta(X^{(j)})^2}
            {\prod_{j=1}^{p-1} \Delta(X^{(j)}, X^{(j+1)})}
         \end{align}
where $w^{(j)}$ is a certain weight function on $\Gamma_j$ for
$j=1, \ldots, p$, with $w^{(1)} = w_1$.
Here
\[ X^{(j)} = (x_1^{(j)}, x_2^{(j)}, \ldots, x_{N_j}^{(j)}) \in \Gamma_j^{N_j} \]
and $\tilde{Z}_n$ is a normalizing constant.

\section{Weak asymptotics}

An important question about a sequence of polynomials with
increasing degrees, is about the asymptotic behavior as the degree
tends to $\infty$.

\subsection{Vector equilibrium problems}
To describe the weak asymptotics of the type II MOPs, as well
as the convergence for the Hermite-Pad\'e rational approximation problems,
vector equilibrium problems were identified that are
relevant for the Angelesco and Nikishin systems. Here we show that
these equilibrium problems have a natural interpretation in
terms of the joint p.d.f.'s of the Angelesco and Nikishin ensembles.

We assume that we are considering MOPs $P_{\vec{n}}$ for a sequence of multi-indices $\vec{n}$
such that $n = |\vec{n}| \to \infty$ and $n_j \to \infty$ for
every $j=1, \ldots, p$ in such a way that
\begin{equation} \label{eq:raylimit}
    \frac{n_j}{n}  \to r_j \qquad \text{for } j = 1, \ldots, p.
    \end{equation}
The limiting ratios $r_j$ should satisfy
\[ 0 < r_j < 1, \qquad \sum_{j=1}^p r_j. \]

Here and in the following we use
\[ I(\mu) = \iint \log \frac{1}{|x-y|} d\mu(x) d\mu(y) \]
to denote the logarithmic energy of the measure $\mu$, and
\[ I(\mu,\nu) = \iint \log \frac{1}{|x-y|} d\mu(x) d\nu(y),
\]
which is the mutual logarithmic energy of the two measures $\mu$ and $\nu$.

For a discrete measure we introduce the reduced logartithmic energy
\[ I^*(\mu) = \iint\limits_{x \neq y} \log \frac{1}{|x-y|} d\mu(x) d\mu(y) \]
Note that, if
\[ \nu_X = \sum_{k=1}^n \delta_{x_k} \]
is the point counting measure of $X = (x_1, \ldots, x_n) \in \mathbb R^n$,
then
\[ \log \Delta(X)^2 = -  I^*(\nu_X). \]

\subsection{Angelesco system}

The Angelesco ensemble has the joint p.d.f.\ \eqref{eq:Angelescopdf}.
A configuration of points $X$ in an Angelesco ensemble
is of the form
\begin{equation} \label{eq:AngelescoX}
    X = (X^{(1)}, \ldots, X^{(p)}), \qquad \text{where}
    \qquad X^{(i)} = (x^{(i)}_1, \ldots, x^{(i)}_{n_i}) \in \Gamma_i^{n_i}.
    \end{equation}
The most likely configuration minimizes
\begin{equation} \label{eq:AngelescoEnergy1}
    - \sum_{j=1}^p \log \Delta(X^{(j)})^2
    -  \sum_{i=1}^{p-1} \sum_{j=i+1}^p \log \Delta(X^{(i)}, X^{(j)})
    +  \sum_{j=1}^p \sum_{k=1}^{n_j} Q_j(x_k^{(j)})
    \end{equation}
where $w_j = e^{-Q_j}$, among all $X$ of the form \eqref{eq:AngelescoX}.
Introducing the normalized point counting measures
\begin{equation} \label{eq:Angelesconuj}
    \nu_{j} = \frac{1}{n} \nu_{X^{(j)}}, \qquad j=1, \ldots, p,
    \end{equation}
we can rewrite \eqref{eq:AngelescoEnergy1}, after dividing by $n^2$, as
\begin{equation} \label{eq:AngelescoEnergy2}
    \sum_{j=1}^p I^*(\nu_j) +
    \sum_{i=1}^{p-1} \sum_{j=i+1}^p I(\nu_i, \nu_j)
    + \frac{1}{n} \sum_{j=1}^p \int Q_j(x) d\nu_j(x).
\end{equation}

In the limit \eqref{eq:raylimit} we forget about the
discreteness of the measures $\nu_j$. Then instead of
minimizing \eqref{eq:AngelescoEnergy2}
among all vectors of measures $(\nu_1, \ldots, \nu_p)$
with $\nu_j$ a measure of total mass $n_j/n$ on $\Gamma_j$
of the form \eqref{eq:Angelesconuj},
we come to minimize the energy functional
\begin{equation} \label{eq:AngelescoEnergy}
    E(\mu_1, \mu_2, \ldots, \mu_p) = \sum_{j=1}^p I(\mu_j) +
    \sum_{j=1}^{p-1} \sum_{k=j+1}^p I(\mu_j, \mu_k)
    \end{equation}
    among all vectors of measures $(\mu_1, \ldots, \mu_p)$ with
\begin{equation} \label{eq:AngelescoNormalization}
    \supp(\mu_j) \subset \Gamma_j \qquad \text{and} \qquad
    \int d\mu_j = r_j.
    \end{equation}

Under the assumption that each $\Gamma_j$ is compact,
and that $w_j(x) \geq 0$ almost everywhere on $\Gamma_j$,
Gonchar and Rakhmanov \cite{GR1} showed that the zeros of
the type II MOP $P_{\vec{n}}$ are distributed according
to the minimizer $(\mu_1, \mu_2, \ldots, \mu_r)$ for the
vector equilibrium problem. More precisely,  for every $j=1, \ldots, p$,
there are
$n_j$ simple zeros of $P_{\vec{n}}$ in $\Gamma_j$, say
$x_1^{(j)}, \ldots, x_{n_j}^{(j)}$, and the
normalized zero counting measure
\[ \nu_j = \frac{1}{n} \sum_{k=1}^{n_j} \delta_{x_k^{(j)}} \]
converges in the limit \eqref{eq:raylimit} weakly to $\mu_j$.

Under the same conditions, it seems likely that the vector of
normalized counting measures $(\nu_1, \ldots, \nu_p)$
(as in \eqref{eq:Angelesconuj}) of a random point
\eqref{eq:AngelescoX} from an Angelesco ensemble
tends to the vector of nonrandom measures $(\mu_1, \ldots, \mu_p)$,
almost surely, but this has not been established rigorously.

In a situation of varying weights in an Angelesco ensemble,
such as for example
\[ w_j(x) =  e^{-n V_j(x)}, \qquad x \in \Gamma_j, \quad j=1, \ldots, p, \]
we have to add  external field terms to \eqref{eq:AngelescoEnergy}
and the relevant energy functional becomes
\begin{equation} \label{eq:AngelescoEnergyFields}
    \sum_{j=1}^p I(\mu_j) +
    \sum_{j=1}^{p-1} \sum_{k=j+1}^p I(\mu_j, \mu_k)
    + \sum_{j=1}^p \int V_j(x) d\mu_j(x).
    \end{equation}
This concept is well known in the orthogonal polynomial case,
see \cite{ST} for the standard reference on logarithmic potential theory
with external fields.

\subsection{Nikishin system}

Similar considerations apply to the joint p.d.f.\
\eqref{eq:NikishinExtendedPDF} in the extended
Nikishin ensemble.
Here the most likely configuration
\[ X = (X^{(1)}, X^{(2)}, \ldots, X^{(p)}) \]
where
\[ X^{(j)} = (x_1^{(j)}, x_2^{(j)}, \ldots, x_{N_j}^{(j)}) \in \Gamma_j^{N_j} \]
is the one that minimizes
\begin{equation} \label{eq:NikishinEnergy1}
    -   \sum_{j=1}^p \log \Delta(X^{(j)})^2
    +  \sum_{j=1}^{p-1} \log \Delta(X^{(j)}, X^{(j+1)})
    +   \sum_{j=1}^p \sum_{k=1}^{N_j}  Q_j(x_k^{(j)})
    \end{equation}
where $w^{(j)} = e^{-Q_j}$.
In terms of the normalized point counting measures
\begin{equation} \label{eq:Nikishinuj}
    \nu_{j} = \frac{1}{n} \nu_{X^{(j)}}, \qquad j=1, \ldots, p,
    \end{equation}
the expression \eqref{eq:NikishinEnergy1} is, after dividing by $n^2$,
\begin{equation} \label{eq:NikishinEnergy2}
    \sum_{j=1}^p I^*(\nu_j) -
    \sum_{j=1}^{p-1} I(\nu_j, \nu_{j+1})
    + \frac{1}{n} \sum_{j=1}^p \int Q_j(x) d\nu_j(x).
\end{equation}

The measure $\nu_j$ has total mass
\[ \int d\nu_j = \frac{N_j}{n} = \frac{1}{n} \sum_{i=j}^p n_i, \]
and in the limit \eqref{eq:raylimit} we have
\[ \frac{N_j}{n} \to \sum_{i=j}^p r_i. \]

So in the limit \eqref{eq:raylimit} we are led
to the following energy functional
\begin{equation} \label{eq:NikishinEnergy}
    E(\mu_1, \mu_2, \ldots, \mu_p) =
    \sum_{j=1}^p I(\mu_j) -  \sum_{j=1}^{p-1} I(\mu_j, \mu_{j+1})
    \end{equation}
among all vectors of measures $(\mu_1, \ldots, \mu_p)$ with
\begin{equation} \label{eq:NikishinNormalization}
    \supp(\mu_j) \subset \Gamma_j \qquad \text{and} \qquad
    \int d\mu_j = \sum_{i=j}^p r_i.
\end{equation}
In the special case that
$n_1 = n_2 = \cdots = n_p$ we have that all $r_i$ are equal to $1/p$,
and then the normalizations are
\[ \int d\mu_j = 1 - \frac{j-1}{p}, \qquad j = 1, \ldots, p. \]

\subsection{Interaction matrix}

In both an Angelesco system and a Nikishin system
we minimize
    \[ E(\mu_1, \mu_2, \ldots, \mu_p)
        = \sum_{j=1}^p \sum_{k=1}^p c_{jk} I(\mu_j, \mu_k) \]
with a certain positive definite interaction matrix $C = (c_{jk})$.

The Angelesco interaction matrix is a full matrix
\begin{equation} \label{eq:AngelescoInteraction}
    C = \begin{pmatrix}
        1 & \frac{1}{2} & \frac{1}{2} & \cdots & \frac{1}{2} & \frac{1}{2} \\
        \frac{1}{2} & 1 & \frac{1}{2} & \cdots & \frac{1}{2} & \frac{1}{2} \\
        \frac{1}{2} & \frac{1}{2} & 1 &      &  & \vdots \\
        \vdots & \vdots & & \ddots & & \vdots \\
        \frac{1}{2} & \frac{1}{2} & \cdots & \cdots & 1 &  \frac{1}{2} \\
        \frac{1}{2} & \frac{1}{2} & \cdots & \cdots & \frac{1}{2} & 1
        \end{pmatrix}
        \end{equation}
and the Nikishin interaction matrix is a tridiagonal matrix
\begin{equation} \label{eq:NikishinInteraction}
    C = \begin{pmatrix} 1 & -\frac{1}{2} & 0 & \cdots & \cdots  & 0 \\
        -\frac{1}{2} & 1 & -\frac{1}{2} &  &   & \vdots   \\
        0 & - \frac{1}{2} & 1 & \ddots    & & \vdots \\
        \vdots    &   & \ddots & \ddots & -\frac{1}{2} & 0 \\
        \vdots &  &  & -\frac{1}{2} & 1 &  -\frac{1}{2} \\
        0 & \cdots  & \cdots & 0 & -\frac{1}{2} & 1
        \end{pmatrix}.
        \end{equation}

The Angelesco interaction matrix \eqref{eq:AngelescoInteraction}
shows the repulsion that exists between $\mu_j$ and $\mu_k$
when $j \neq k$. However, the repulsion is only half as strong
as the inner repulsion between each measure $\mu_j$ itself.

Since the non-zero off-diagonal entries in \eqref{eq:NikishinInteraction}
are negative, there is an attraction in a Nikishin ensemble
between two consecutive measures $\mu_j$ and $\mu_{j+1}$ for $j=1, \ldots, p-1$.
The measures $\mu_j$ and $\mu_k$ with $|j-k| \geq 2$ do not
interact.
The Nikishin interaction also arises in the asymptotic analysis of
eigenvalues of banded Toeplitz matrices \cite{DuK1} as well as in the two-matrix model \cite{DuK2}
where it appears with both an external field and an upper constraint.

\bibliographystyle{amsalpha}

\end{document}